\documentclass[11pt,longbibliography]{article}

\usepackage[right=3cm,left=3cm,textheight=22cm]{geometry}

\usepackage{amsmath}
\usepackage{amsthm}
\usepackage{amssymb}
\usepackage{amsfonts}
\usepackage{paralist}
\usepackage{url}
\usepackage{color}
\usepackage{graphicx}
\usepackage{subfig}
\usepackage{tikz}

\allowdisplaybreaks

\usepackage[colorlinks=true,linkcolor=blue,citecolor=blue,urlcolor=blue,breaklinks]{hyperref}
\hypersetup{pdftitle={Existence of Compactly Supported Global Minimisers for the Interaction Energy}}
\hypersetup{pdfauthor={José A. Cañizo, José A. Carrillo \& Francesco Saverio Patacchini}}

\numberwithin{equation}{section}

\newtheorem{thm}{Theorem}[section]

\newtheorem{proposition}[thm]{Proposition}

\theoremstyle{definition}

\theoremstyle{remark}

\newcommand{\R}{\mathbb{R}}
\newcommand{\re}{\mathbb{R}}
\newcommand{\real}{\mathbb{R}}
\newcommand{\ren}{\mathbb{R}^N}

\newcommand{\prob}{\mathcal{P}}
\newcommand{\B}{\mathcal{B}}
\def\supp{\mbox{supp}}

\def\eps{\varepsilon}
\newcommand{\dive}{\mathrm{div}}
\def\na{\nabla}
\def\qed{\,\unskip\kern 6pt \penalty 500
\raise -2pt\hbox{\vrule \vbox to8pt{\hrule width 6pt
\vfill\hrule}\vrule}\par}

\newcommand{\be}{\begin{equation}}
\newcommand{\ee}{\end{equation}}
\newcommand{\bfig}{\begin{figure}}
\newcommand{\efig}{\end{figure}}
\newcommand{\bt}{\begin{table}}
\newcommand{\et}{\end{table}}
\newcommand{\bc}{\begin{center}}
\newcommand{\ec}{\end{center}}
\newcommand{\ba}{\begin{array}}
\newcommand{\ea}{\end{array}}
\newcommand{\bes}{\begin{equation*}}
\newcommand{\ees}{\end{equation*}}

\begin{document}

\title{ Some Free Boundary Problems involving\\ Nonlocal Diffusion and Aggregation}

\author{
Jos\'e Antonio Carrillo\footnote{Department of Mathematics, Imperial College London, London SW7 2AZ.}, Juan Luis V\'{a}zquez\footnote{Departmento de Matem\'aticas, Universidad Aut\'onoma de Madrid, 28049-Madrid, Spain.}}

\maketitle


\begin{abstract}
We report on recent progress in the study of evolution processes involving degenerate parabolic equations what may exhibit free boundaries. The equations we have selected follow to recent trends in diffusion theory: considering anomalous diffusion with  long-range effects, which leads to fractional operators or other operators involving kernels with large tails; and the combination of diffusion and aggregation effects, leading to delicate long-term equilibria whose description is still incipient.
\end{abstract}

\section{Introduction}\label{sec.intro}

The systematic mathematical study of free boundary problems is quite recent compared to the present importance of the field, and goes back to the last half of the past century. In this period it has developed in many directions, which combine modeling of quite diverse natural phenomena, particularly in diffusion, wave mechanics and elasticity. A main topic in this broad landscape has always been the study of evolution processes where the state equations governing the different phases are nonlinear parabolic equations, possibly degenerate or singular at certain points or for certain values of the unknowns. This is the case of the well-known models like Stefan Problem for the common evolution of two immiscible fluids, the Porous Medium Equation, and other models to be mentioned below.

Free boundary models combine the difficulty of the analysis of the nonlinear PDEs with the difficulty in locating the separating  interface between the phases,  or free boundary in the absence of a theory of classical solutions that is generally unavailable.


So in principle we can solve the problem in some generalized sense but neither the state variables are guaranteed to be continuous or the free boundaries to be reasonable hyper-surfaces.  Even when the initial and boundary data are smooth, nonlinear processes involving free boundaries may develop singularities in finite time both in the solution and the free boundary. This is therefore a difficult problem for both analysis and geometry. In the last 50 years a considerable progress has been achieved and this is reported in different monographs, and in its recent trends in other articles of this volume. Several new models have been actively pursued in recent years with various levels of mathematical rigor versus practical applicability. In this article we will discuss free boundary problems arising in two new scenarios, nonlocal diffusion and  aggregation/swarming processes.

\section{Nonlocal diffusion}

The classical theory of diffusion is expressed mathematically by means of the heat equation, and more generally by parabolic equations, normally of the linear type; such approach  has  had an enormous success and is now a foundation stone in science and technology. The last half of the past century has witnessed intense activity and progress in the theories of nonlinear diffusion, examples being the Stefan Problem, the Porous Medium Equation, the $p$-Laplacian equation, the  Total Variation Flow, evolution problems of Hele-Shaw type, the Keller-Segel chemotaxis system, and many others. Mean curvature flows also belong to this broad category. Reaction diffusion has also attracted considerable attention.

In the last decade there has been a surge of activity focused on the use of so-called fractional diffusion operators to replace the standard Laplace operator (and other kinds of elliptic operators with variable coefficients), with the aim of further extending the theory by taking into account the presence of the long range interactions that occur in a number applications. The new operators do not act by pointwise differentiation but by a global integration with respect to a singular kernel;  in that way the {\sl nonlocal character} of the process is expressed.

More generally, research on non-local operators has a tradition related to stochastic processes, and has witnessed a rapid expansion in the last decade when it has attracted the attention of PDE experts who brought new problems and techniques into the field. The area poses new challenges to both pure and applied mathematicians, and is now at the interface of at least three wide fields: functional analysis, stochastic processes and partial differential equations, while providing a new paradigm in scientific modeling.  This interaction between mathematics and applications is giving rise to new concepts and methods, and is expected to produce new challenging mathematical problems for many years to come.


\

\noindent {\bf \large Fractional operators}

\noindent Though there is a wide class of interesting nonlocal operators under scrutiny, both in theory and applications, a substantial part of the current work deals with diffusion modeled by the so-called fractional Laplacians. We recall that the fractional Laplacian operator is a kind of isotropic differentiation operator of order $2s$, for some $s\in (0,1)$, that can be conveniently defined  through its Fourier Transform symbol, which is $|\xi|^{2s}$.  Thus, if $g$ is a function in the Schwartz class in $\ren$, $N\ge 1$, we write  $(-\Delta)^{s}  g=h$ if
\begin{equation*}
\widehat{h}\,(\xi)=|\xi|^{2s} \,\widehat{g}(\xi)\,,
\end{equation*}
so that for $s=1$ we recover the standard Laplacian. This definition allows for a wider range of parameters $s$. The interval of interest for fractional diffusion is $0<s\le 1$, and for $s< 1$  we can also use the integral representation
\begin{equation}\label{def-riesz}
(-\Delta)^{s}  g(x)= C_{N,2s }\mbox{
P.V.}\int_{\mathbb{R}^N} \frac{g(x)-g(z)}{|x-z|^{N+{2s}
}}\,dz,
\end{equation}
where P.V. stands for principal value and $C_{N,\sigma}$ is a normalization constant, with precise value that is found in the literature. In the limits $s\to 0$ and $s\to 1$ it is possible to recover respectively the identity or the standard minus Laplacian, $-\Delta$, cf. \cite{Brezis}. Remarkably, the latter one cannot be represented by a nonlocal formula of the type \eqref{def-riesz}.
It is also useful to recall that the operators $(-\Delta)^{-s}$, $0< s<1$,  inverse of the former ones,  are  given by standard convolution expressions:
\begin{equation*}
(-\Delta)^{-s}  g(x)= C_{N,-2s }\int_{\mathbb{R}^N}
\frac{g(z)}{|x-z|^{N-{2s} }}\,dz,
\end{equation*}
in terms of the usual Riesz potentials. Basic references for these operators are the books by
Landkof \cite{Landkof} and Stein \cite{Stein}. A word of caution: in the literature we often find the notation $\sigma=2s$, and then the desired interval is $0<\sigma<2$. According to that practice, we will sometimes use $\sigma$ instead of $s$. These and other definitions  are equivalent when dealing with the Laplacian on the whole space $\re^N$.
The correct definitions for the operators defined on bounded domains admit several options that are being investigated at the moment, \cite{BSV, SV1}. The interest in these fractional operators has a long history in Probability and other applied sciences. The systematic study of the corresponding PDE models with fractional operators is relatively recent, and many of the results have been established in the last decade, see for instance the presentations \cite{DPV11,Kass1,Valdinoc,VazAbel,VazDCDS}, where further references can be found.

A part of the current research concerns linear or quasilinear equations of elliptic type. This is a huge subject with well-known classical references. Even if it includes Free Boundary Problems, mainly of the Obstacle Type, see \cite{CaffSilv07}, it  will not be discussed here in itself for lack of space, because we want to present free boundaries in evolution problems.

\


\noindent {\bf \large Linear evolution processes of anomalous type}

\noindent  We will study evolution models that arise as variants of the heat equation paradigm, including however long-range effects. Thus, the difference between the standard and the fractional Laplacian consists in taking into account long-range interactions (jumps) instead of the usual interaction driven by close neighbors. The change of  model explains characteristic new features of great importance, like enhanced propagation with the appearance of fat tails at long distances: such tails great differ from the typical exponentially small tails of the standard diffusion, and even more from the compactly supported solutions of porous medium flows. Moreover, the space scale of the propagation of the distribution is not proportional to $t^{1/2}$ as in the Brownian motion, but to another power of time, that can be adjusted in the model; this is known as {\sl anomalous diffusion.}  Anomalous diffusion is nowadays intensively studied, both theoretically and experimentally since it explains a number of phenomena in several areas of physics, finance, biology, ecology, geophysics, and many others, which can be  briefly summarized as having non-Brownian scaling, see e.g. \cite{Applebaum, ContTankov2004, Woy2001}.
The fractional Laplacian operators of the form $(-\Delta)^{\sigma/2}$, $\sigma\in(0,2)$, are actually the infinitesimal generators of stable L\'{e}vy processes \cite{Applebaum, Bertoin, CKS2010}.
The standard linear evolution equation involving fractional diffusion is
\begin{equation*}
\frac{\partial u}{\partial t}+(-\Delta)^{s}(u)=0\,,
\end{equation*}
which is  the main model for anomalous diffusion.  The equation is solved with the aid of well-known tools, like Fourier transform. Posed in the whole space, it  generates a semigroup of ordered contractions in $L^1(\mathbb{R}^N)$ and has the integral representation
\begin{equation*}
u(x,t)=\int_{\mathbb{R}^N}K_s(x-z,t)f(z)\,dz\,,
\end{equation*}
where $K_s$ has Fourier transform  $\widehat K_s(\xi,t)=e^{-|\xi|^{2s} t}$. This means that, for $0<s<1$, the kernel $K_s$ has the self-similar form
$ K_s(x,t)=t^{-N/2s}F(|x|/t^{1/2s}) $
 for some profile function $F=F_s$ that is positive and decreasing, and it behaves at infinity
like $F(r)\sim r^{-(N+2s)}$, \cite{Blumenthal-Getoor}. When $s=1/2$, $F$ is explicit:
\begin{equation*}
F_{1/2}(r)={C}\,(a^2+r^2)^{-(N+1)/2}\,.
\end{equation*}
If $s=1$ the function $K_{s=1}$ is the Gaussian heat kernel, which has a negative square exponential tail, i.e., a completely different long-distance behavior.

From the point of view of the present investigation, this fractional linear model and other variants, though useful and heavily studied, fail to meet the requirement that we are looking for, that is, generating free boundaries in the evolution. On the contrary, the solutions have large densities at long distances.

\section{Nonlinear and nonlocal diffusion models}

A main feature of current research in the area of PDEs  is the interest in nonlinear equations and systems. There are a number of models of  evolution equations that  have been proposed as nonlinear counterparts of the linear fractional heat equation, and combine Laplace operators and nonlinearities in different ways. Let us mention the two popular in the recent PDE literature, and how they relate to our stated goal.

\smallskip

\noindent $\bullet$ {\bf Model I.} A quite natural option is to consider the equation
\begin{equation}\label{fpme}
\partial_t u +(-\Delta)^{s}(u^m)=0
\end{equation}
with $0<s<1$ and $m>0$. This is mathematically the simplest fractional version of the standard Porous Medium Equation (PME) $\partial_t u =\Delta(u^m)\,,$ that is recovered as the limit $s\to 1$ and has been extensively studied, cf. \cite{ArBk86, Vapme}. We will call equation \eqref{fpme} the {\sl Fractional Porous Medium Equation, FPME,} as proposed in our first works on the subject \cite{pqrv1, pqrv2}, in collaboration with de A. de Pablo,  F. Quir\'os and A. Rodr\'{\i}guez. From the point of view of our present study, this model is quite interesting because it offers a possible balance between a porous medium nonlinearity, which for $m>1$ implies finite propagation and free boundaries, and a fractional-type diffusion where the speed on propagation is infinite and fat tails arise. The question is: which effect wins? In the papers \cite{pqrv1, pqrv2} we have answered the question in the sense of infinite propagation and no free boundaries (we are considering nonnegative solutions of the evolution problem).

Nowadays, a reasonably complete theory of nonlinear diffusion with long-range effects has been developed for this model, see for instance \cite{pqrv3, VazBar2012, VazVol, BV2012}, or similar models like \cite{CJ11}, but we have to look elsewhere for free boundaries.

\smallskip

\noindent $\bullet$  {\bf Model II.} A second model of nonlinear diffusion with fractional operators has been studied by the author in collaboration with Luis Caffarelli and does give rise to the occurrence of finite propagation and free boundaries. This alternative model is derived in a more classical way from the Porous Medium Equation since it is based on the usual Darcy  law (i.\,e., the velocity of the particles is assumed to be the gradient of a pressure function, ${\bf v}=-\nabla p$) with the novelty that the pressure is related to the density by an inverse fractional Laplacian operator, $p={\cal K}u$. Putting these two facts into the continuity equation, $\partial_t u+\nabla\cdot(u\bf v)=0$, the model takes the form
\begin{equation}\label{PMEfp}
u_t=\nabla (u\, \nabla {\mathcal K} u),
\end{equation}
a  nonlinear fractional diffusion equation of porous medium type where $\mathcal K$ is the Riesz operator that typically expresses the inverse to the fractional Laplacian, ${\mathcal K} u=(-\Delta)^{-s}u$. This has been studied by Caffarelli and the author in \cite{CV1, CVobs}. In dimension $N=1$ this model was studied by Biler, Karch and Monneau \cite{BKM, BIK} as a model for the propagation of dislocations as proposed by Head \cite{Head}. Since it was extensively reported in the survey paper \cite{VazAbel} we will only give point out to the established theory of weak solutions, and the lack of uniqueness results in several dimensions, and concentrate on the existence and properties of the free boundaries.  Self-similar solutions  exist and their existence and properties are very illustrative of the phenomenon of finite propagation.

\subsection{Finite propagation for Model II. Solutions with compact support}

One of the most important features of the porous medium equation
and other related degenerate parabolic equations is the property
of finite propagation, whereby  compactly supported initial data
$u_0(x)$ gives rise to solutions $u(x,t)$ that have the same
property for all positive times, i.e., the support of $u(\cdot,t)$
is contained in a ball $B_{R(t)}(0)$ for all $t>0$.
 One possible proof in the case of the PME is by constructing
explicit weak solutions exhibiting that property (i.e.,
having a free boundary) and then using the comparison principle,
that holds for that equation.
 Since we do not have such a general principle here, we have to
devise a comparison method with a suitable family of  ``true
supersolutions'', which are in fact some quite excessive
supersolutions. The technique has to be adapted to the peculiar
form of the integral kernels involved in operator ${
\cal K}_s$.

We begin with $N=1$ for simplicity. We assume that our solution
$u(x,t)\ge 0$ has bounded initial data $u_0(x)=u(x,t_0)\le M$ with
compact support and is such that {\sl
 $u_0$ \quad is below the parabola \quad $a(x-b)^2$, $a,b>0$,
with graphs strictly separated.} We may assume that $u_0$ is
located under the left branch of the parabola. We take as
comparison function
$$
\sl  U(x,t)=a(Ct-(x-b))^2,
$$
which is a traveling wave moving to the right with speed $C$ that will be taken big enough.
Then we argue at the first point and time where $u(x,t)$ touches  the left branch of the parabola
$U$ from below. The key point is that if $C$ is large enough such contact cannot exist. The formal idea
is to write the equation as
$$
u_t= u_xp_x+ up_{xx}
$$
and observe that at the contact we have $u_t\ge U_t=2aC(Ct-x+b)$, while $u_x=U_x=-2a(Ct-x+b)$, so the first
can be made much bigger than the second by increasing $C$. The influence of $p_x$ and $p_{xx}$ as well as $u$ is controlled, and then we conclude that the equation cannot hold if $C$ is large enough. The argument can be translated for several dimensions.
Here are the  detailed results proved in \cite{CV1}.

\begin{thm}\label{thmex1} Let $0<s<1/2$ and assume that $u$ is a bounded solution of equation
{\rm (\ref{PMEfp})} with $0\le u(x,t)\le L$, and  $u_0$ lies  below a  function of  the form
$U_0(x)=Ae^{-a|x|}, \ A,a>0.$ If $A$ is large then there is a constant $C>0$ that depends only on $(N,s,a, L, A)$ such that for any $T>0$ we will have the comparison
\begin{equation*}
u(x,t)\le Ae^{Ct-a|x|} \quad \mbox{for all $x\in \mathbb{R}^N$ and all $0<t\le T$.}
\end{equation*}
If  $1/2\le s<1$ a similar statement is true but $C$ is not a constant but some increasing function of time.
\end{thm}

\noindent A similar finite propagation result is true in several space dimensions. The study of the corresponding free boundaries is an open topic. Let us mention in passing that much effort has been devoted in papers \cite{CV1,CSV,CV3} to prove the basic sequence of regularity results, according to which initial data in $L^1$ give rise to solutions that belong to $L^p$, $1\le p\le \infty$ for all positive time. Moreover,
 there exists a positive constant $C$ such that for every $t>0$
\begin{equation*}
\sup_{x\in\mathbb{R}^N}|u(x,t)|\le C(N,s)\,t^{-\alpha }\|u_0\|_{L^1(\mathbb{R}^N)}^{\gamma}
\end{equation*}
with $\alpha=N/(N+2-2s)$, $\gamma=(2-2s)/((N+2-2s)$. Finally, it is proved that bounded weak solutions  $u\ge 0 $ of the initial value problem are uniformly continuous on bounded sets of $s<1$. Indeed, they are $C^\alpha$ continuous  with a uniform modulus.

\subsection{Self-similar solutions for Model II}

Next we want to  study of the large time behaviour of this model following paper \cite{CVobs}, since it shows some light on the way finite propagation works. The first step is constructing the self-similar solutions that will serve as attractors. Surprising result: their density will be compactly supported, while their pressure will be positive everywhere with the typical fat tails far away.

\smallskip

\noindent  {\bf Rescaling for the FPME.}
 Inspired by the asymptotics of the standard porous medium equation,
we define the rescaled (also called renormalized) flow through the
transformation
\begin{equation}\label{scaling}
\sl  u(x,t)=(t+1)^{-\alpha} v(x/(t+1)^{\beta}, \tau)\normalcolor
\end{equation}
with new time $\tau=\log(1+t)$. We also put $y=x/(t+1)^{\beta}$ as
rescaled space variable. In order to cancel the factors including
$t$ explicitly, we  get the condition on the exponents
$ \alpha + (2-2s)\beta =1\,. $ Here we use the homogeneity of ${\cal K}$ in the form
$
\sl  ({\cal K} u)(x,t) = t^{-\alpha+2s\beta} ({\cal
K}v)(y,\tau).
$
From physical considerations we also impose the law that states
conservation of (finite) mass, which amounts to the condition
$\alpha= N \beta$, and  In this way  we arrive at
 the precise value for the exponents:
\begin{equation*}
\sl  \beta=1/(N+ 2-2s), \quad
\alpha=N/(N+ 2-2s).\normalcolor
\end{equation*}
 We also arrive at the  {\sl nonlinear, nonlocal Fokker-Planck equation}
\begin{equation}\label{ren.eq}
v_\tau=\nabla_y\cdot(v\,(\nabla_y {\cal K}(v)+\beta y))
\end{equation}
which is the equation for the {\sl renormalized flow.} In all the above calculations the
factor $(t+1)$ can be replaced by $t+t_0$ for any $t_0>0$, or even by plain $t$.

\medskip

\noindent  {\bf Stationary renormalized solutions}. It is important to concentrate on the stationary states of
the new equation, i.\,e., on the solutions  $V(y)$ of
\begin{equation*}
\sl  \nabla_y\cdot(V\,\nabla_y (P+ a |y|^2))=0, \quad \mbox{with \ } \ P={\cal
K}(V).\normalcolor
\end{equation*}
where $a=\beta/2$, and $\beta$ is defined just above. Since we are
looking for asymptotic profiles of the standard solutions
of the FPME we also want $V\ge 0$ and integrable. The simplest
possibility is integrating once to get
\begin{equation}\label{baren.eq2}
\sl   V\,.\nabla_y (P+  a |y|^2))=0, \quad P={\cal K}(V), \quad
V\ge 0.\normalcolor
\end{equation}
The first equation gives an alternative choice that reminds us of the
complementary formulation of the obstacle problems.
  Indeed, if we solve the {\bf Obstacle problem with
fractional Laplacian} we will obtain a unique solution $P(y)$ of
the problem:
\begin{equation}\label{obst}
\begin{array}{l}
P\ge f, \quad V=(-\Delta)^{s} P\ge 0; \\
\mbox{either } \ P=f \ \mbox{or } \ V=0.
\end{array}
\end{equation}
with $0<s<1$.  In order for solutions of \eqref{obst} to be also solutions of \eqref{baren.eq2} \normalcolor  we have to choose as obstacle $f(y)=C- a \,|y|^2,$
where $C$ is any positive constant and $a=\beta/2$. Note that $-\Delta f=2Na=\alpha$. For
uniqueness we also need the condition $P\to 0$ as $|y|\to\infty$.
Fortunately, the corresponding theory had been developed by Caffarelli and collaborators,  cf. \cite{CSS, Silv07}. The solution is unique and belongs to the space $H^{-s}$ with
pressure in $H^{s}$. Moreover, it is shown that the  solutions have $P\in C^{1,s}$ and $V\in
C^{1-s}$.

Note that for $C\le 0$ the solution is trivial, $P=0$, $V=0$, hence we choose $C>0$. We also note the pressure is defined but for a constant, so that we could maybe take as pressure $\widehat P=
P-C$ instead of $P$ so that $\widehat P=0$; but this does not simplify things since $P\to 0$ implies that $\widehat P\to -C$ as $|y|\to\infty$. Keeping thus the original proposal, we get a one parameter family of
stationary profiles that we denote $V_C(y)$. These solutions of the obstacle problem
produce correct weak solutions of the fractional PME equation with
initial data a multiple of the Dirac delta for the density, in the form
\begin{equation*}
U_C(x,t)=t^{-\alpha}V_C(|x|t^{-\beta}).
\end{equation*}
It is what we can call the
source-type or Barenblatt solution for this problem, which is a
profile $V\ge 0$. It is positive in the {\sl  contact set } of the
obstacle problem, which has the form ${\cal C}=\{|y|\le R(C)\}$,
and is zero outside, hence it has compact support.

On the other hand, the rescaled pressure $P(|y|)$ is always
positive and decays to zero as $|y|\to\infty$ according to
fractional potential theory, cf. Stein \cite{Stein}.
The rate of decay of $P$ as $|y|\to\infty$ turns out to be $P=O(|y|^{2s-N})$, a fat tail.

\medskip

\noindent{\bf Exact calculation of density profiles.}  Biler, Karch and Monneau \cite{BKM} studied the existence and stability of self-similar solutions in one space dimension. Recently,  Biler, Imbert and Karch \cite{BIK} obtain the explicit formula for a multi-dimensional self-similar solution in the form
\begin{equation*}
U(x,t)=c_1t^{-\alpha} (1-x^2t^{-2\alpha/N})_+^{1-s}
\end{equation*}
with $\alpha=N/(N+2-2s)$ as before.  The derivation uses an important identity for fractional Laplacians which is found in Getoor \cite{Get61}: \ $(-\Delta)^{\sigma/2}(1-y^2)_+^{\sigma/2}=K_{\sigma,N}$ \ if $\sigma\in (0,2]$. Here we must take $\sigma=2(1-s)$. According to our previous calculations $\Delta P=-\alpha$ on the coincidence set, hence  $c_1=\alpha/K_{\sigma,N}$. Let us work a bit more: using a similar scaling to \eqref{scaling}, we arrive at the following one-parameter family of self-similar solutions
\begin{equation*}
U(x,t;C_1)=t^{-\alpha} (C_1-k_1\,x^2t^{-2\alpha/N})_+^{1-s}
\end{equation*}
where $k_1=c_1^{1/(1-s)}$ and $C_1>0$ is a free parameter that can be fixed in terms of the mass of the solution $M=\int U(x,t;C_1)\,dx$. This is the family of densities that corresponds to the pressures obtained above as solutions of the obstacle problem. Let us finally mention that explicit formulas for self-similar solutions to related fractional equations have also recently been obtained in \cite{Huang}.

\subsection{Asymptotic behaviour for Model II}

\noindent The next step is  to prove that these profiles are attractors for the rescaled flow.
The studay is done in terms of the rescaled flow that is much more workable than the original evolution.  The following result is proved in \cite{CVobs}.

\begin{thm}\label{thm.asymp} Let $u(x,t)\ge 0$ be a weak solution of the Cauchy Problem for equation
\eqref{PMEfp} with bounded and integrable initial data such that $u_0 \ge 0$ has finite entropy, i.e.,
$$
\int_{\mathbb{R}^N} (v\,{\cal K}(v) + \beta |y|^2v)\,dy < \infty\,.
$$
Let $v(y,\tau)$ be the corresponding rescaled solution to \eqref{ren.eq}. As $\tau\to\infty$ we have
 \begin{equation*}
 v(\cdot,\tau)\to V_C(y) \quad \mbox{in \  $ L^1(\mathbb{R}^N)$ and also in $L^\infty(\mathbb{R}^N)$}.
 \end{equation*}
 The constant $C$ is determined by the rule of mass equality: $\int_{\mathbb{R}^N} v(y,\tau)\,dy=\int _{\mathbb{R}^N} V_C(y)\,dy$. In terms of function $u$, this translates into the convergence as $t\to \infty$
 \begin{equation*}
 u(x,t)-U_C(x,t)\to 0 \quad \mbox{in } L^1(\mathbb{R}^N), \quad t^{\alpha}|u(x,t)-U_C(x,t)|\to 0 \quad \mbox{uniformly in $x$}\,.
 \end{equation*}
 \end{thm}
\noindent The constants $\alpha$ and $\beta$ are the self-similar exponents defined before. The proof uses entropy dissipation methods that adapt perfectly to the problem. Sharper results on this convergence are given in \cite{CHSV} where more references can be found. See also below in this respect.

\smallskip

\noindent {\bf Extension of the nonlinear fractional models.} {\bf 1)} There is way of generalizing the two previous Models {\bf  I} and {\bf II}, to accept two exponents $m$ and $p$:
$$
\partial_t u +\nabla(u^{m-1}\nabla (-\Delta)^{-s}u^{p})=0\,,
$$
so that the comparison of both models happens on symmetric terms. The question we want to solve is deciding between finite and infinite propagation in terms of the exponents $m$ and $p$. Progress on this issue is reported in \cite{StanTesoVazquez, StanTesoVazquez2}.

\noindent {\bf 2) Limits. } The limit of Model {\bf II} when $s\to1$ is quite interesting since one  obtains a variant of the equation for the evolution of vortices in superconductivity studied in \cite{CRS}, E \cite{WE, Liz, AmSr}. The understanding of this limit has been done in collaboration with Sylvia Serfaty \cite{SerVaz}, and is related to work by Bertozzi et al. on aggregation models \cite{BCL, BLL}. This is a good point to connect to the next issue. On the other hand, the limit $m\to\infty$ of Model 1 leads to a free boundary problem of the Mesa Type reported in \cite{VazBar2012, VazMesa2014}.

\section{Aggregation and swarming}

The collective behavior of individuals in animal grouping, also called swarming, has recently been modelled by mean field PDEs, see \cite{KCBFL} and the references therein. They consist in continuum descriptions of Individual Based Models (IBMs), see \cite{Ao,Reyn,HW} for applications to fish schools and software simulation of swarms. The motion of individuals in IBMS is governed by systems of ODEs imposing an asymptotic speed, cruise speed, for particles. One of the most celebrated models of this kind was introduced in \cite{Dorsogna}. The authors propose to introduce a social force between individuals based on the empirical observation that there are two basic mechanisms of interaction: a inner repulsion zone to keep a comfort area or to avoid collisions and an outer attraction area since individuals do want to socialize and be close to each other. These two tendencies are modeled via a pairwise potential directed along the position between the particles and depending on the interparticle distance. Therefore, this model consists in
\begin{flalign}\label{micro}
\begin{split}
   \frac{d x_i}{d t} &= v_i \,, \\
   \frac{d v_i}{d t} &= \alpha v_i  - \beta v_i \vert v_i\vert^2
   - \sum_{j \neq i} \nabla W(x_i-x_j)\,,
\end{split}
\end{flalign}
where $x_i,v_i\in\R^N, i=1,\dots,M$ are the positions and velocities of the individual particles and $\alpha, \beta$ are effective values for self-propulsion and friction forces, see \cite{Dorsogna} for more discussions. These models lead to kinetic equations in the limit of large number of particles as derived in \cite{CDP}. Some of the patterns that the simulation of these IBMs produce are the so-called {\sl flock patterns}. The flock solution corresponds to the case in which all particles move with the same velocity vector $v_i=v_0$ with $|v_0|^2=\alpha/\beta$ while the positions $\hat x_i\in\R^N, i=1,\dots,M$ satisfy that the forces are balanced, that is,
$$
\sum_{j \neq i} \nabla W(\hat x_i-\hat x_j) =0\,.
$$
The flocking solution corresponds then to the translational movement of the shape formed by the positions  $\hat x_i\in\R^N, i=1,\dots,M$ with velocity vector $v_0$. When the number of particles tends to infinity, these discrete flocking solutions are approximating continuum density profiles $\hat \rho$ satisfying $\nabla W \ast \hat \rho =0$ at least over the support of $\hat \rho$. In fact, the connection between the second-order microscopic model \eqref{micro} and the first-order model given by
\begin{equation} \label{aggregation}
\frac{d x_i}{d t} = - \sum_{j \neq i} \nabla W(x_i-x_j)\,,
\end{equation}
has been clarified in \cite{Albietal,CHM}. There, it is shown that the stability of flock solutions for \eqref{micro} is equivalent to the stability of the spatial shape of the flock as steady solution to \eqref{aggregation}. A similar result at the level of the continuum equations is lacking. However, it makes sense to find all the spatial configurations leading to flocks or equivalently all stationary solutions to \eqref{aggregation}. Among those, the most stable configurations are of special interest, in view of the connection through the stability of the solutions at least at the discrete level. The system \eqref{aggregation} has a natural energy that is dissipated along its evolution given by
$$
E[x_1,\dots,x_M]=\frac12 \sum_{i,j: i \neq j} W(x_i-x_j)\,.
$$
The continuum version of this energy, suitably scaled, leads to
\begin{equation}\label{eq:contE}
E[\mu] = \frac12\int_{\mathbb{R}^N\times\mathbb{R}^N} W(x-y)
d\mu(x)d\mu(y)\,
\end{equation}
where $\mu$ is a probability measure. This energy is dissipated along the evolution of the continuum mean-field limit of the system of ODEs in \eqref{aggregation}, called the aggregation equation
\begin{equation}
  \label{eq:aggregation}
  \partial_t \rho = \mbox{div} (\rho (\nabla W * \rho))\,,
\end{equation}
where $\rho$ is the density function of the measure $\mu$.
Therefore, minimizers of the energy \eqref{eq:contE} on the set of probability measures should be among the most stable configurations. Finally, let us comment that both the nonlocal diffusion equation \eqref{PMEfp} and its rescaled version \eqref{ren.eq} are particular instances of the aggregation equation \eqref{eq:aggregation} for suitably chosen potentials $W$. We will see in the rest of this section a summary of the results concerning the minimizers of the interaction energy \eqref{eq:contE}, specially in its relation to free boundary problems, particularly to the obstacle problem.

\medskip

\noindent {\bf Connection to optimal transport.}
We start by noting that the discrete system \eqref{aggregation} has the structure of a gradient flow with respect to the Euclidean distance. In fact, this structure is shared by the continuum equation \eqref{eq:aggregation} in a suitable sense. We need to be precise about the definition of distance between probability measures. We recall that $\prob(\real^N)$ is  the set of Borel probability measures on $\real^N$ and we denote by  $\B(\real^N)$ the family of Borel subsets of $\real^N$. The support of a measure $\mu \in \prob(\real^N)$ is the closed set defined by
\begin{equation*} 
\supp (\mu):= \{x \in \real^N: \mu(B(x,\epsilon))>0  \text{ for all } \epsilon>0 \}\,.
\end{equation*}
A family of distances between probability measures has been classically
introduced by means of optimal transport theory, we will review briefly some of
these concepts, we refer to \cite{V} for further details. A probability
measure $\pi$ on the product space $\real^N \times \real^N$ is said to be a
transference plan between $\mu \in \prob(\real^N)$ and $\nu \in \prob(\real^N)$ if
\begin{equation}\label{marginal}
\pi(A \times \real^N)=\mu(A) \quad \text{and} \quad \pi(\real^N
\times A)=\nu(A)
\end{equation}
for all $A \in \B(\real^N)$. If $\mu, \nu \in \prob(\real^N)$,
then
$$
\Pi(\mu,\nu):=\{ \pi \in \prob(\real^N \times \real^N):
\eqref{marginal} \text{ holds for all }A \in \B(\real^N)\}
$$
denotes the set of admissible transference plans between $\mu$ and
$\nu$. Informally, if $\pi \in \Pi(\mu,\nu)$ then $d \pi(x,y)$
measures the amount of mass transferred from location $x$ to
location $y$. We recall that the Euclidean Wasserstein distance $d_2$ between two
measures $\mu$ and $\nu$ is defined by
\[
    d_2^2(\mu,\nu)=\inf_{\pi\in\Pi(\mu,\nu)}\left\lbrace \iint_{\R^N\times \R^N}
    |x-y|^2 d\pi(x,y)\right\rbrace.
\]
Note that $d_2(\mu,\nu)<\infty$ for $\mu,\nu\in \prob_2(\real^N)$
the set of probability measures with finite moments of order $2$. It is classical that the equation \eqref{eq:aggregation} is a gradient flow of the energy \eqref{eq:contE} in the space of measures $\prob_2(\real^N)$ endowed with the distance $d_2$, see \cite{Carrillo-McCann-Villani03,Carrillo-McCann-Villani06}. This property is essential to obtain many qualitative properties of the solutions and to deal with possible singular interaction potentials, see \cite{BCL,CDFLS,CDFLS2,BT2,BLL,BCLR,BCLR2,BCY,CCH2} and the references therein.

In our context of minimizing the energy functional \eqref{eq:contE}, this means that we look for local minimizers in $d_2$. (Local) Minimizers of \eqref{eq:contE} should correspond to equilibrium configurations for the evolution equation obtained by steepest descent of the energy. However, being a functional on probability measures, the steepest descent has to be understood in the Wasserstein sense as in \cite{Otto} by writing \eqref{eq:aggregation} as
\begin{equation}\label{eq:main}
\frac{\partial \mu_t}{\partial t} = \dive \left[\left(\nabla
\frac{\delta E}{\delta\mu}\right)\mu_t\right] = \dive \left( \mu_t \nabla\psi_t
\right) \qquad x\in\R^N \, , t>0.
\end{equation}
This evolution equation can make sense in the set of probability measures depending on the regularity of the potential $W$. However, if the potential is singular at the origin the well-posedness only happens in $L^p$ spaces, see \cite{BLR,CDFLS,CCHau} for more details.
For steady configurations, we expect $\na \psi = 0$ on the support of $\mu$, due
to the formal energy dissipation identity for solutions, i.e.,
$$
\frac{d}{dt} E[\mu_t] = -\int_{\R^N} |\nabla \psi_t|^2\,d\mu_t
$$
where $\mu_t$ is any solution at time $t$ of \eqref{eq:main} and $\psi_t=W\ast
\mu_t$ its associated potential. Therefore, the points on the support of a local
minimizer $\mu$ of the energy $E$ should correspond to critical points of its
associated potential $\psi$. This fact is made rigorous in:

\begin{proposition}[{\rm \cite[Theorem 4]{BCLR2}}]\label{prop:minsupp}
Assume that $W$ is a non-negative lower semi-continuous function in
$L^1_{loc}(\R^N)$. If $\mu$ is a $d_2$-local minimizer of the energy, then
the potential $\psi$ satisfies the Euler-Lagrange conditions given by

(i) $\psi(x)=(W\ast \mu)(x)=2E[\mu]$ $\mu$-a.e.

(ii) $\psi(x)=(W\ast \mu)(x)\leq 2E[\mu]$ for all $x\in
\supp(\mu)$.

(iii) $\psi(x)=(W\ast\mu)(x)\geq 2E[\mu]$ for a.e. $x\in\R^N$.
\end{proposition}

\noindent These conditions simplify to
\begin{align}
\psi(x)=(W\ast \mu)(x) &= 2E[\mu]\quad \text{for a.e.}\quad x\in\supp(\mu)\,,\nonumber\\
\psi(x)=(W\ast \mu)(x) &\geq 2E[\mu]\quad \text{for a.e.} \quad x\in
\R^N\setminus \supp(\mu)\,,\label{cucu}
\end{align}
if $\mu$ is absolutely continuous with respect to the Lebesgue measure. These Euler-Lagrange conditions can also be interpreted in terms of game theory, see \cite{Cardia,BC,DLR}. In fact, they imply that the potential $\psi$ achieves its global minima on the support of $\mu$, and therefore $\mu$ is also characterized as
$$
\int_{\mathbb{R}^N} \psi(x) d \mu(x)= \min_{\nu\in\prob_2(\real^N)} \int_{\mathbb{R}^N} \psi(x) d \nu(x)\,.
$$
In game theory the probability measure $\mu$ is interpreted as the strategy chosen by the agent $x$, the potential is interpreted as encoding the interactions between the agents, and the minimization is thought about the optimal condition (the game) that the agents are looking for. Therefore, the Euler-Lagrange conditions are interpreted in game theory as the lack of interest of any player in the game to change of strategy leading to the concept of Nash equilibria. As a summary, local minimizers of the potential energy are among Nash equilibria within the game theoretical viewpoint.

\begin{figure}[htb]
\begin{center}
\begin{tikzpicture}
\draw[black,->] (-0.15,0) -- (8,0) node[anchor=north] {$a$};
\draw[thick,dashed, color=green!80] (1.2,1.2) -- (8,1.2) node[black,anchor=west] {$\,\,2-N$};
\draw[black,dashed] (1.2,2) -- (8,2) node[anchor=west] {$\,\,3-N$};
\draw[black,dashed] (1.2,3) node[anchor=east] {$2$} -- (3,3);
\draw[black]	 (1.2,0) node[anchor=north] {$2-N$}
	              (3,0) node[anchor=north] {2}
	              (8,0) node[anchor=west] {$\,\,-N$}
	              (3,0) node[anchor=north] {2}
                    (2.1,0) node[anchor=north] {0};
\fill[color=black!15] (8,1.2) rectangle (2.1,0);
\draw[black,->] (1.2,-0.15) -- (1.2,4) node[anchor=east] {$b$};
\draw[thick,red,dashed] (3,0) -- (3,1.2);
\draw[very thick,blue,dashed] (3,1.2) -- (3,2);
\draw[black,dashed] (3,2) -- (3,3);

\filldraw[black] (3,1.2) circle (2pt);

\draw[black] (0,0) -- (4,4) node[anchor=west] {$a=b$};

\end{tikzpicture}
\end{center}

\caption{Parameters for minimizers of the interaction energy with potential of the form \eqref{eq:pot0} given by densities, see the text for a full explanation.} \label{fig:param}
\end{figure}
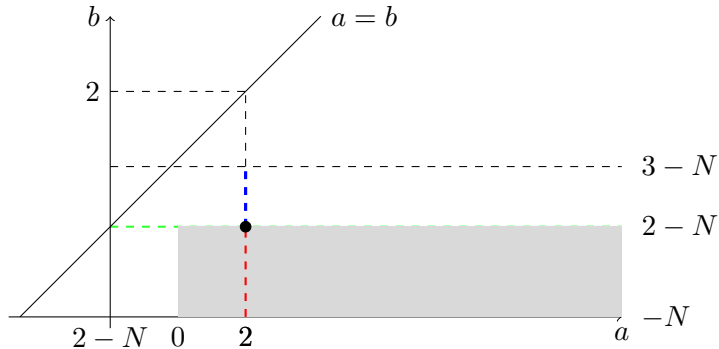

Let us concentrate on the quest of local minimizers for repulsive/attractive potentials of the form
\begin{equation}\label{eq:pot0}
W (x)=\frac{|x|^{a}}{a} - \frac{|x|^b}{b}\,,
\end{equation}
with $a>b>-N$ to have a repulsive at the origin and attractive at infinity potential which is locally integrable. We allow $a$ or $b$ to be zero with the understanding that
${|x|^{0}}/{0}\sim \log|x|$.

As mentioned before, the case $a=2$ and $b=2-N$, black point in Figure \ref{fig:param}, is of particular interest. It corresponds to Newtonian repulsion and it has been repetitively rediscovered that the unique, up to translations, global minimizer of the interaction energy
is the characteristic function of a Euclidean ball with suitable radius \cite{Fr,PH,CGZ}. The uniqueness, up to translations, of the global minimizer for more singular than Newtonian
repulsion, $a=2$, $b=2s-N$ with $0<s<1$ (red dashed line in Figure \ref{fig:param}), was obtained in \cite{CVobs} via the connection to a classical obstacle problem, see also \cite{CGZ}. This strategy was also used in \cite{SerVaz} to treat again the case $s=1$ for the evolution problem as in \cite{BLL}. If the parameters are in the range of the green line in Figure \ref{fig:param}, i.e. $a>2$ or $2-N<a<2$ with $b=2-N$, the existence and uniqueness of compactly supported radial minimizers of the interaction energy was obtained in \cite{FHK,FH}. Summing up, the solution support is limited by a free boundary, whose location (a sphere) is to be determined as part of the problem.

Moreover, the main result in \cite{CDM} shows that in the range $a>0$ with $b=2s-N$ with $0<s\leq1$, grey area in Figure \ref{fig:param}, local minimizers in $d_2$ of the interaction energy \eqref{eq:contE} are absolutely continuous with respect to the Lebesgue measure, and their density function lies in  $L^\infty_{loc}(\R^N)$ when $s=1$ and in $C^\alpha_{loc}(\R^N)$ when $s\in(0,1)$. These results are obtained again by exploiting the connection between the Euler-Lagrange conditions for local minimizers and classical obstacle problems. They can be generalized for repulsive at the origin potentials behaving like $|x|^b$ with $b=2s-N$, $0<s\leq 1$, and under quite general conditions of the behavior of the potential $W$ outside the origin. In all these cases, we know that the local minimizers are compactly supported and determined by a free boundary whose regularity is still unknown in full generality.

Let us conclude this part by mentioning that whenever we can fully characterize the unique, up to translations, global minimizers of \eqref{eq:contE}, they are the long-time asymptotics for the aggregation equation \eqref{eq:aggregation}, see \cite{BLL,CV1,CVobs,SerVaz} for the known different cases (and Section 3 above). These cases are essentially reduced to $a=2$ and $b=2s-N$ with $0<s\leq1$, and even decay rates are known in one dimension \cite{CFP,CHSV}. Finally, let us remark that finding sufficient conditions on the interaction potential $W$ to ensure the existence of compactly supported global minimizers is a very interesting question that has been recently solved in \cite{CCP}, see also \cite{SST,CFT,CCH}. However, we do not know yet about their properties in terms of asymptotic stability for the aggregation equation  \eqref{eq:aggregation}.

\medskip

\noindent {\bf More on the connection  to free boundary problems and the obstacle problem}.
 We will now focus on a topic already introduced in Section 3 (b) above. A classical obstacle problem \cite{CSS} for a differential operator $\mathcal{L}$ consists in finding a function $\Phi:\R^N \to \R$ such that
\begin{equation}
  \label{eq:obstacle}
  \left\{
    \ba{l}
      (\Phi - f)(\mathcal{L} \Phi - g) = 0
      \\
      \Phi \geq f
      \\
      \mathcal{L} \Phi \geq g
    \ea\right.
  \text{on $\R^N$,}
\end{equation}
where $f,g: \R^N \to \R$ are just given functions. Here, $f$ is the obstacle to $\Phi$, and thus whenever $f$ does not touch the obstacle, then the first condition implies that $\Phi$ should be a solution to the PDE: $\mathcal{L} \Phi= g$. To make the connection to the Euler-Lagrange condition clearer, let us take the particular case of the classical potential encountered in
semiconductors \cite{CF,CDMS} given by
$$
W(x)=c_N \frac{|x|^{2-N}}{2-N} + \frac{|x|^2}{2}
$$
with $N\geq 1$. Here, $c_N$ is the constant such that the first part of the potential is the fundamental solution of the Laplacian operator $\Delta$ in $N\geq 1$. Assume that
$\rho \in L^1\cap L^\infty (\R^N)$ is a $d_2$-local minimizer of the energy associated to this potential, then the potential associated to $\rho$, given by $\psi = W\ast \rho$ is well defined as a function that satisfies in the sense of distributions
$$
\Delta \psi = \Delta W \ast \rho = \rho + N\,,
$$
and thus, $\Delta \psi \geq N$ in $\R^N$ and $\Delta \psi = N$ outside the support of $\rho$. In conclusion, the potential $\psi$ associated to the local minimizer satisfies a classical obstacle problem in $\R^N$ with operator $\mathcal{L}=\Delta$, $g=N$, and $f=2E[\rho]$ due to \eqref{cucu}. Let us also remark that the obstacle problem \eqref{obst} associated to the scaled version of the fractional diffusion equation \eqref{ren.eq} can also be recast in the form \eqref{eq:obstacle}. We will do this next but in more generality.

The basic idea for the Laplacian can be generalized and used to get properties on the potentials associated to local minimizers of a larger family of interaction potentials $W$ without resorting to the local minimizers themselves. Let us review quickly the main strategy in \cite{CDM}. The key fact is that if the repulsive strength at zero of the potential is given by ${|x|^b}/{b}$, $b=2s-N$ with $0<s\leq1$, then one can show that the potential $\psi = W\ast \mu$ associated to any $d_2$-local minimizer $\mu$ is a continuous function, see \cite{CDM}. This allows us to deduce that for any point $x_0\in\supp(\mu)$
\begin{equation}\label{eq:psimine}
 \psi(x)\geq \psi(x_0) \mbox{ for all } x\in B_\eps(x_0)
\end{equation}
holds. Here, $\eps$ is the size of the ball in $d_2$ distance where $\mu$ is a local minimizer. Moreover, for the points in the support, we also get
\begin{equation}\label{eq:psiminequal}
\psi (x) = \psi(x_0) \quad  \mbox{ in } B_\eps(x_0) \cap\supp(\mu).
\end{equation}
Since $\frac{|x|^b}{b}$ in this range is the fundamental solution, up to constants, of the fractional diffusion operator ${\mathcal L}=-(-\Delta)^{s}$, we can write in the distributional sense that
$$
(-\Delta)^{s} \psi = \mu +(-\Delta)^{s} \frac{|x|^{a}}{a} * \mu \quad \mbox{ in } \mathcal D'(\R^N)\,.
$$
In particular, since $\mu$ is a non-negative measure, we deduce
$$
(-\Delta)^{s} \psi \geq (-\Delta)^{s} \frac{|x|^{a}}{a} * \mu \quad \mbox{ in $B_\eps(x_0)$}.
$$
Furthermore, if $x\in B_\eps(x_0)$ is such that $\psi(x) >\psi(x_0)$,
\eqref{eq:psiminequal} implies that $x\notin \supp(\mu)$. We deduce
$$
(-\Delta)^{s} \psi = (-\Delta)^{s} \frac{|x|^{a}}{a} * \mu \quad \mbox{ in } \mathcal D'(\R^N) \mbox{ in } B_\eps(x_0)\cap \{\psi>\psi(x_0)\}\,.
$$
Collecting \eqref{eq:psimine}-\eqref{eq:psiminequal}, we have, at least formally, the following result:

\begin{thm}[{\rm \cite{CDM}}]\label{prop:obstacle}
For all $x_0\in\supp (\mu)$, the potential function $\psi$ associated to a $d_2$-local minimizer of the energy \eqref{eq:contE} is equal, in $B_\eps(x_0)$, to the unique solution of the obstacle problem
\begin{equation}\label{eq:obstacle2l}
\left\{
\begin{array}{rll}
\varphi &\geq C_0,\quad & \mbox{ in } B_\eps(x_0) \\
(-\Delta)^{s} \varphi &\geq  - F(x), \quad & \mbox{ in } B_\eps(x_0) \\
(-\Delta)^{s} \varphi &=  - F(x), & \mbox{ in } B_\eps(x_0)\cap\{\varphi >C_0\} \\
\varphi &= \psi, & \mbox{ on  } \partial B_\eps(x_0),
\end{array}
\right.
\end{equation}
where $C_0=\psi(x_0)$ and $F(x)=-(-\Delta)^{s} \frac{|x|^{a}}{a}*\mu$. Furthermore,
we can recover the local minimizer of the energy $\mu$ as $\mu = -\Delta \psi + F$.
\end{thm}
Most of the results mentioned above follow this strategy to make use of the regularity and existence theorems for obstacle problems available in the literature, see \cite{C,C1,Blank,Silv07,CSS}. Let us remark that in the case of a power law potential \eqref{eq:pot0} with $a=2$ and $0<s\leq 1$, $F(x)$ is constant. The range corresponding to $0<s<1$ and $a=2$ corresponds to the obstacle problem treated by Caffarelli and V\'azquez in \cite{CVobs}. As before, Theorem \ref{prop:obstacle} is true for more general potentials behaving as the repulsive power $|x|^b$ at the origin with suitable growth conditions at infinity, see \cite{CDM}.

It is an open question how to use obstacle problem techniques in cases in which the repulsion is not too singular, for instance for power-law potentials when \eqref{eq:pot0} with $a=2$ and $b=2s-N$ with now $s>1$. The main problem we face is that the inverse operator of the convolution is no longer as nice as it was for fractional diffusion. However, we are currently investigating the case when $s=1+\epsilon$ with $0<\epsilon<1$, i.e., in the range corresponding to the blue line in Figure \ref{fig:param}. In this particular case, we consider an intermediate obstacle problem by defining $\Psi$ to be the solution to
\begin{equation}\label{obsnew}
\begin{array}{l}
\Psi\ge c_1, \quad V=(-\Delta)^{\epsilon} \Psi\ge 0; \\
\mbox{either } \ \Psi=c_1 \ \mbox{or } \ V=0\,,
\end{array}
\end{equation}
which again has nice regularity properties since $0<\epsilon<1$. Now, we define our potential $\psi$ as the Newtonian potential associated to $\Psi$, that is, we solve $\Psi=-\Delta \psi$. The idea is that by integrating twice the solution of the obstacle problem \eqref{obsnew}, we get a candidate to be a solution to \eqref{eq:obstacle2l}.

\section*{Acknowledgement}
JAC acknowledges support from projects MTM2011-27739-C04-02, 2009-SGR-345 from Ag\`encia de Gesti\'o d'Ajuts Universitaris i de Recerca-Generalitat de Catalunya, the Royal Society through a Wolfson Research Merit Award, and the Engineering and Physical Sciences Research Council (UK) grant number EP/K008404/1. JLV acknowledges support from Spanish Project MTM2011-24696, and a long-term stay at the Isaac Newton Institute, Cambridge, UK.

{\small
\bibliographystyle{royalagain}

}

\end{document}